\documentclass{owrart}

\usepackage{epsfig,psfrag}
\usepackage{xypic}

\def\R{\mathbb R}
\def\N{\mathbb N}
\def\M{\overline{M}_{0,n}(\R)}
\def\Mp{\overline{M}_{0,n+1}(\R)}

\def\Fl{\mathcal F\!\ell_m}
\def\A{\mathcal A}
\def\C{\mathcal C}
\def\acts{\hspace{.1cm}{\setlength{\unitlength}{.30mm}\linethickness{.09mm}
                        \begin{picture}(8,8)(0,0)\qbezier(7,6)(4.5,8.3)(2,7)\qbezier(2,7)(-1.5,4)(2,1)\qbezier(2,1)(4.5,-.3)(7,2)
                                                 \qbezier(7,6)(6.1,7.5)(6.8,9)\qbezier(7,6)(5,6.1)(4.2,4.4)
                        \end{picture}\hspace{.1cm}}}
\def\acted{\hspace{.1cm}{\setlength{\unitlength}{.30mm}\linethickness{.09mm}
                        \begin{picture}(8,8)(0,0)\qbezier(1,6)(3.5,8.3)(6,7)\qbezier(6,7)(9.5,4)(6,1)\qbezier(6,1)(3.5,-.3)(1,2)
                                                 \qbezier(1,6)(1.9,7.5)(1.2,9)\qbezier(1,6)(3,6.1)(3.8,4.4)
                        \end{picture}\hspace{.1cm}}}
\def\lllarrow{  \hspace{.05cm}\mbox{\,\put(0,-3){$\leftarrow$}\put(0,1){$\leftarrow$}\put(0,5){$\leftarrow$}\hspace{.45cm}}}
\def\squareA{{
    \setlength{\unitlength}{.15mm}
    \linethickness{.08mm}
    \begin{picture}(10,10)(0,0)    \qbezier(0,0)(0,0)(10,10)
    \qbezier(0,0)(0,0)(10,0)    \qbezier(0,0)(0,0)(0,10)
    \qbezier(10,0)(10,0)(10,10)    \qbezier(0,10)(0,10)(10,10)
        \end{picture}
    }}
\def\squareB{{
    \setlength{\unitlength}{.15mm}
    \linethickness{.08mm}
    \begin{picture}(10,10)(0,0)    \qbezier(10,0)(10,0)(0,10)
    \qbezier(0,0)(0,0)(10,0)    \qbezier(0,0)(0,0)(0,10)
    \qbezier(10,0)(10,0)(10,10)    \qbezier(0,10)(0,10)(10,10)
    \end{picture}
    }}
\def\pentagonA{{
    \setlength{\unitlength}{.15mm}
    \linethickness{.08mm}
    \begin{picture}(13,15)(-4,-.1)    \qbezier(0,0)(0,0)(9,0)
    \qbezier(0,0)(0,0)(-3,9)    \qbezier(-3,9)(-3,9)(4.5,14.5)
    \qbezier(4.5,14.5)(12,9)(12,9)    \qbezier(9,0)(9,0)(12,9)
    \qbezier(0,0)(0,0)(4.5,14.5)
    \qbezier(9,0)(9,0)(4.5,14.5)
    \end{picture}
    }}
\def\pentagonB{{
    \setlength{\unitlength}{.15mm}
    \linethickness{.08mm}
    \begin{picture}(13,15)(-4,-.1)    \qbezier(0,0)(0,0)(9,0)
    \qbezier(0,0)(0,0)(-3,9)    \qbezier(-3,9)(-3,9)(4.5,14.5)
    \qbezier(4.5,14.5)(12,9)(12,9)    \qbezier(9,0)(9,0)(12,9)
    \qbezier(0,0)(0,0)(12,9)    \qbezier(0,0)(0,0)(4.5,14.5)
    \end{picture}
    }}
\def\pentagonC{{
    \setlength{\unitlength}{.15mm}
    \linethickness{.08mm}
    \begin{picture}(13,15)(-4,-.1)    \qbezier(0,0)(0,0)(9,0)
    \qbezier(0,0)(0,0)(-3,9)    \qbezier(-3,9)(-3,9)(4.5,14.5)
    \qbezier(4.5,14.5)(12,9)(12,9)    \qbezier(9,0)(9,0)(12,9)
    \qbezier(0,0)(0,0)(12,9)
    \qbezier(-3,9)(-3,9)(12,9)
    \end{picture}
    }}
\def\pentagonD{{
    \setlength{\unitlength}{.15mm}
    \linethickness{.08mm}
    \begin{picture}(13,15)(-4,-.1)    \qbezier(0,0)(0,0)(9,0)
    \qbezier(0,0)(0,0)(-3,9)    \qbezier(-3,9)(-3,9)(4.5,14.5)
    \qbezier(4.5,14.5)(12,9)(12,9)    \qbezier(9,0)(9,0)(12,9)
    \qbezier(-3,9)(-3,9)(9,0)    \qbezier(-3,9)(-3,9)(12,9)
    \end{picture}
    }}
\def\pentagonE{{
    \setlength{\unitlength}{.15mm}
    \linethickness{.08mm}
    \begin{picture}(13,15)(-4,-.1)    \qbezier(0,0)(0,0)(9,0)
    \qbezier(0,0)(0,0)(-3,9)    \qbezier(-3,9)(-3,9)(4.5,14.5)
    \qbezier(4.5,14.5)(12,9)(12,9)    \qbezier(9,0)(9,0)(12,9)
    \qbezier(-3,9)(-3,9)(9,0)
    \qbezier(9,0)(9,0)(4.5,14.5)
    \end{picture}
    }}

\newtheorem{theorem}{Theorem}

\begin{document}

\begin{talk}{Andre Henriques}
{An action of the cactus group}
{Andr{\'e} Henriques}

Let $\M$ denote the Deligne-Mumford compactification of the moduli space of real curves of genus zero with $n$ marked points.
Its points are the isomorphism classes of stable real curves of genus zero,
that is, curves obtained by glueing $\R \mathbb{P}^1$'s in a tree-like way, and such that each irreducible component has at least 3 special points.
Let $[\Mp/S_n]$ denote the quotient orbifold of $\Mp$ by the action permuting the first $n$ marked points.
In \cite{HK06}, J. Kamnitzer and the author showed that the {\em cactus group} $J_n:=\pi_1([\Mp/S_n])$ 
acts on tensor powers of Kashiwara crystals in a way similar to how the braid group acts on tensor powers of quantum group representations.

The {\em big cactus group} $J'_n$ is the fundamental group of $[\M/S_n]$.
It fits into a short exact sequence
\(
0\to \pi_1(\M)\to J'_n\to S_n\to 0,
\)
and its elements can be represented by movies, such as the following one:
\[
\psfig{file=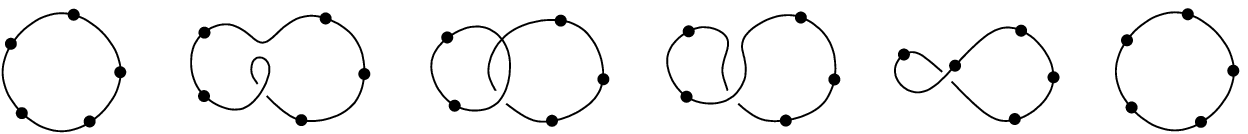, width=8cm}
\]

Let $\Fl:=\big(\put(1,5){$\scriptscriptstyle 1$}\put(4,2.4){$\scriptscriptstyle 1$}
\put(7,-0.2){$\scriptscriptstyle 1$}\put(10,-3.2){$\scriptscriptstyle 1$}\put(9,3.5){$\scriptstyle *$}\hspace{.47cm}\big)\big\backslash SL_m$ 
be the variety of flags $0\subset V_1\subset \cdots \subset V_{m-1}\!\subset \R^m$, equipped with volume forms $\omega_i\in\Lambda^i V_i$.
The goal of this note is to construct an action of $J'_n$ on the totally positive part $\A(n)_{>0}$ of the variety $\A(n):=(\Fl)^n\!/SL_m$.
The space $\A(n)_{>0}$ is a certain connected component of the locus $\A(n)_{\mathit{reg}}\subset \A(n)$, where the flags are in generic position.
One gets similar actions on $\big((N\backslash G)^n\!/G\big){}_{>0}$ for other reductive groups $G$.

The space $\A(n)_{>0}$ was introduced by Fock and Goncharov \cite{FG06}. For $m=2$, it agrees with the Teichm\"uller space of decorated ideal $n$-gons,
that is, the space of isometry classes of hyperbolic $n$-gons with geodesic sides, vertices at infinity, and horocycles around each vertex.
It is also an example of a cluster variety, i.e. it comes with special sets of coordinate systems,
whose transition functions are given by {\em cluster exchange relations} \cite{FZ02}.
For $m=2$, the coordinates are due to Penner \cite{P87}. To each pair $i,j$ of vertices of the $n$-gon,
he associates the quantity $\Delta_{ij}:=\exp(\frac{1}{2}d_{ij})$, 
where $d_{ij}$ denotes the hyperbolic length between the intersection points of the horocycles around $i$ and $j$, and the geodesic from $i$ to $j$.
These coordinates are then subject to the following exchange relations \cite{P87}:
{\psfrag{a}{$d_{ij}$}
\psfrag{b}{$d_{jk}$}
\psfrag{c}{$d_{k\ell}$}
\psfrag{d}{$d_{i\ell}$}
\psfrag{e}{$d_{ik}$}
\psfrag{e'}{$d_{j\ell}$}
\psfrag{EEEEEEEEEEEEEEEEEEEEEEEEEEEEEEEE}{$\displaystyle \Delta_{j\ell}=\frac{\Delta_{ij}\Delta_{k\ell}+\Delta_{jk}\Delta_{i\ell}}{\Delta_{ik}}$.}
\begin{equation}\label{exchange}
\begin{matrix}
\psfig{file=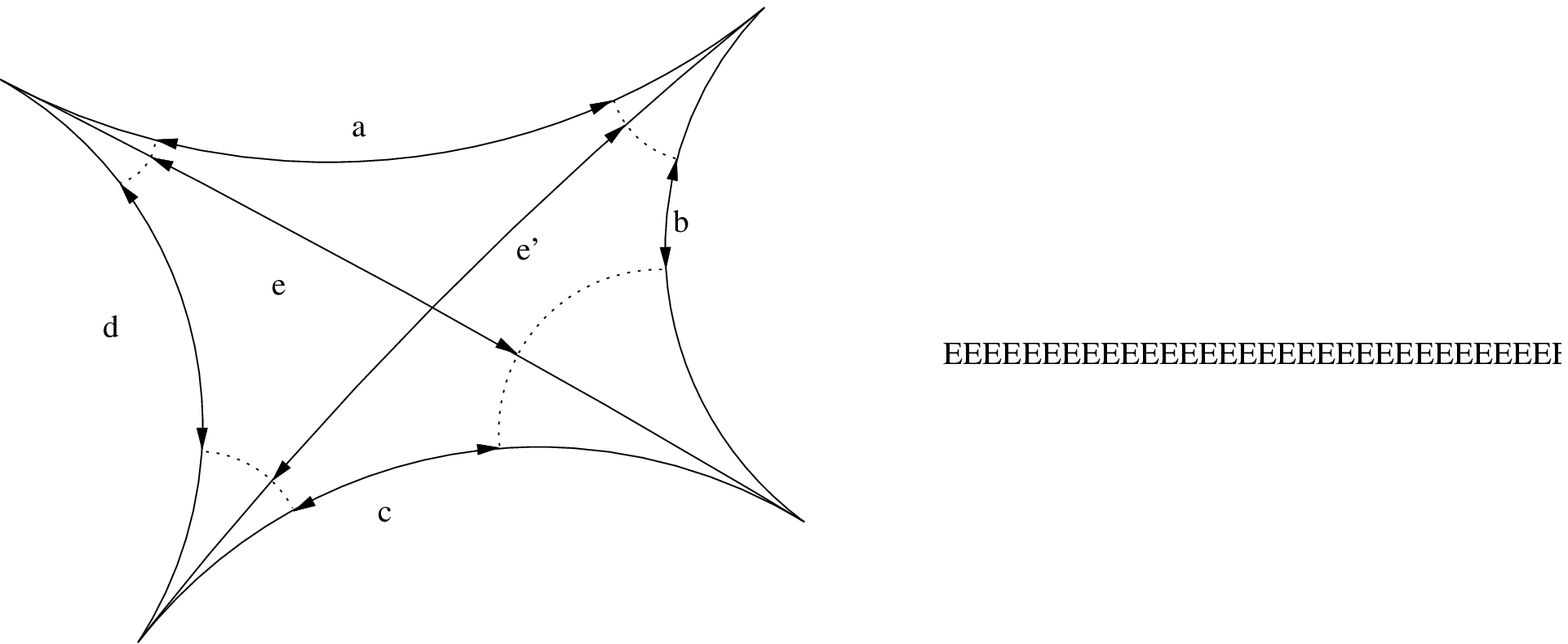, width=6cm}
\end{matrix}
\end{equation}}
\vspace{-.2cm}

\noindent 
For general $m$, the coordinates on $\A(n)$ are indexed by tuples $(i_1,\ldots,i_n)\in\N^n$
whose sum equals $m$, and such that at least two entries are non-zero.
The coordinate $\Delta_{i_1\ldots i_n}$ then assigns to $((V^1_\bullet,\omega^1_\bullet),\ldots,(V^n_\bullet,\omega^n_\bullet))\in(\Fl)^n$ the ratio of 
$\omega^1_{i_1}\wedge\cdots\wedge\omega^n_{i_n}$ with the standard volume form on $\R^m$.
These coordinates satisfy
\[
\begin{split}
\Delta_{\ldots i \ldots j \ldots k \ldots \ell \ldots}
=&\big(
\Delta_{\ldots i+1 \ldots j \ldots k \ldots \ell-1 \ldots}
\cdot\Delta_{\ldots i \ldots j-1 \ldots k+1 \ldots \ell \ldots}
\,+\\
&\Delta_{\ldots i \ldots j \ldots k+1 \ldots \ell-1 \ldots}
\cdot\Delta_{\ldots i+1 \ldots j-1 \ldots k \ldots \ell \ldots}
\big)\big/
\Delta_{\ldots i+1 \ldots j-1 \ldots k+1 \ldots \ell-1 \ldots}\,\,,
\end{split}
\]
which generalizes (\ref{exchange}).
Let $\A(n)_{>0}$ be the locus where all the $\Delta$'s are $>0$.
It is a space isomorphic to $\R_{>0}^{(n-2)\cdot\binom{m+1}{2}+(m+1)-n}$, and each triangulation of the $n$-gon provides such an isomorphism \cite{FG06}.
More precisely, the isomorphism corresponding to a triangulation is given by the coordinates $\Delta_{0\ldots 0i0\ldots 0j0\ldots 0k0\ldots 0}$, 
where $i$, $j$, $k$ are located at the vertices of the triangles.
For example, for $n=8$, $m=4$, and the triangulation \put(4,-3){\psfig{file=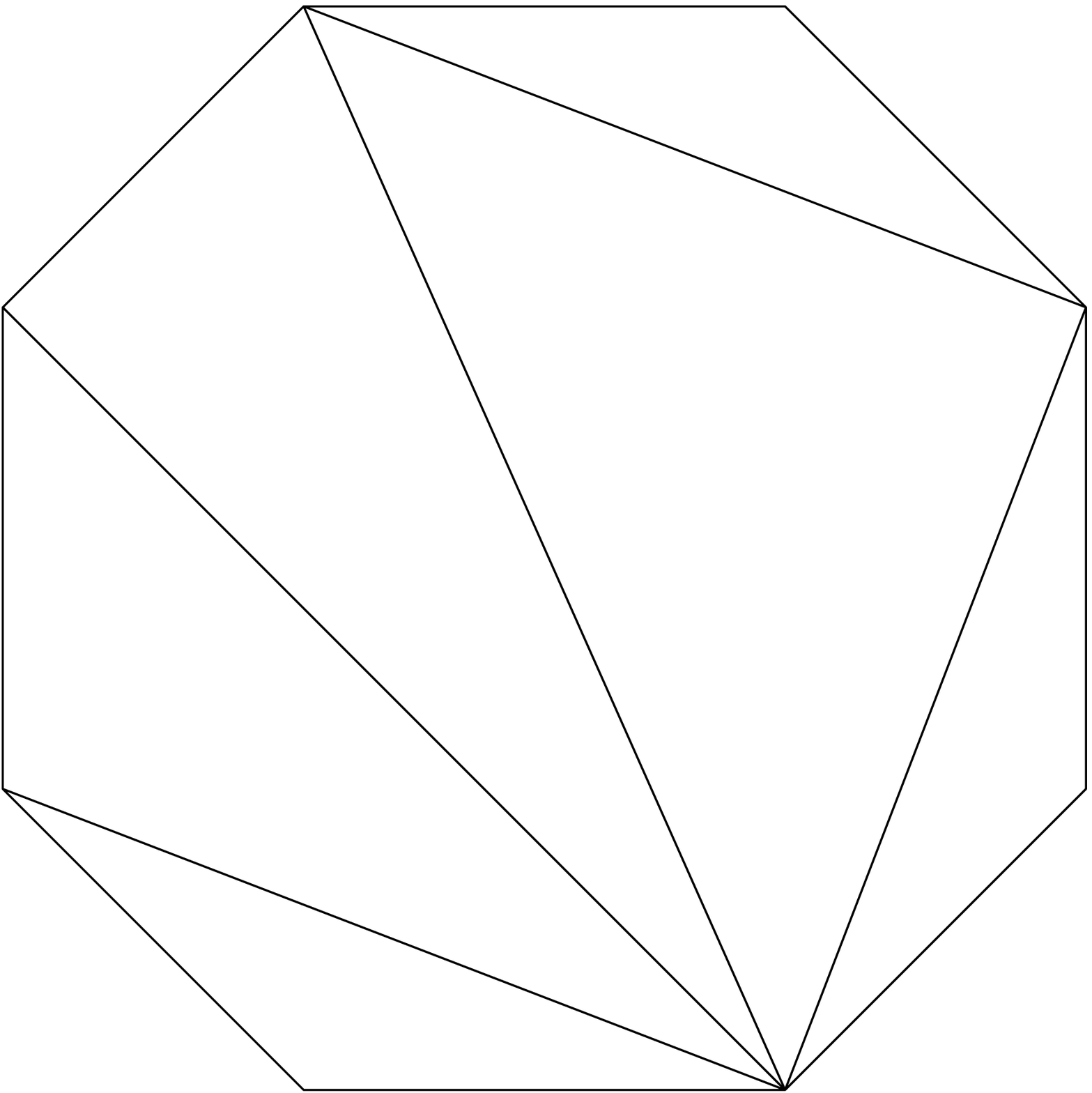, width=0.5cm}}\hspace{.57cm}\phantom{.} of the 8-gon,
the corresponding coordinates on $\A(n)_{>0}$ are in natural bijection with the bullets in the following figure:
\begin{equation}\label{bullets}
\begin{matrix}\psfig{file=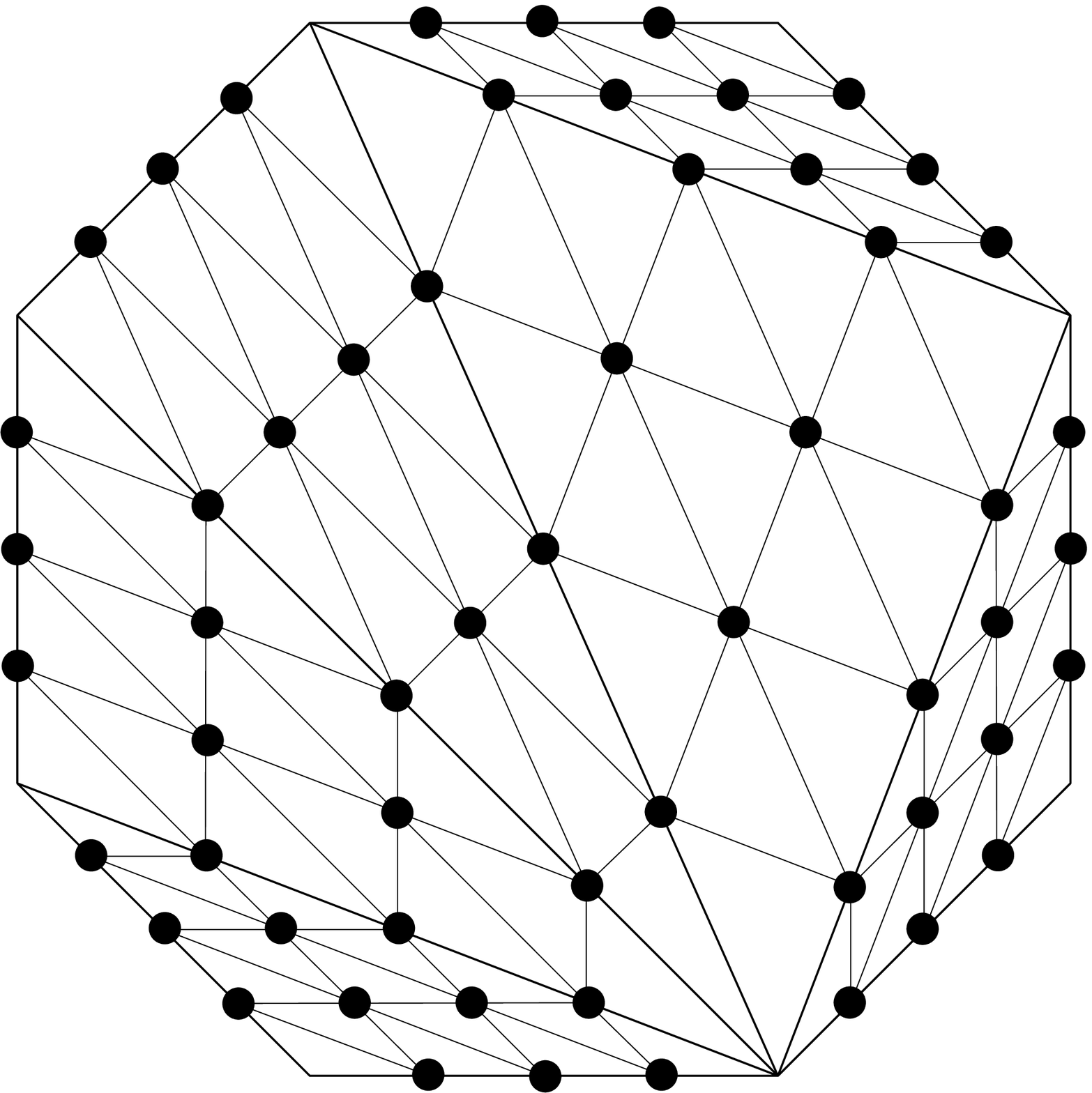, width=2cm}\end{matrix}
\end{equation}

We now explain a general machine for producing actions of $J'_n$ on various spaces. 
Suppose that we are given two manifolds $X_\triangle$ and $X_I$, equipped with maps
\begin{equation}\label{diagram}
r\acts X_\triangle
\put(12,9.5){$\scriptstyle d_1$}
\put(12,1.5){$\scriptstyle d_2$}
\put(5.5,-6){$\scriptstyle d_3=d_0$}
\hspace{.05cm}\mbox{\,\put(0,-3){$-\!\!\!-\!\!\!\longrightarrow$}\put(0,1){$-\!\!\!-\!\!\!\longrightarrow$}\put(0,5){$-\!\!\!-\!\!\!\longrightarrow$}\hspace{1cm}}
X_I\acted\iota
\end{equation}
subject to the relations $r^3=1$, $\iota^2=1$, and $d_i\circ r=r\circ d_{i-1}$.
Such data can then be reinterpreted as a contravariant functor
\(
X_\bullet:\C\rightarrow\{\text{manifolds}\}
\)
from the category $\C:=\{\acts\triangle\lllarrow I\acted\}$, whose two objects 
are the oriented triangle ``$\triangle$'' and the unoriented interval ``$I$'', and whose morphisms are the obvious embeddings and automorphisms.
Let $\widehat\C$ be the category whose objects are the 2-dimensional finite simplicial complexes with oriented 2-faces and connected links, 
and whose morphisms are the embeddings.
There is an obvious inclusion $\C\hookrightarrow \widehat\C$, and every object of $\widehat\C$ can be written essentially uniquely as the
colimit of a diagram in $\C$.
Assuming $d_1\times d_2\times d_3:X_\triangle\to X_I^3$ is a submersion, 
then there is a unique extension of $X_\bullet$ to $\widehat \C$ sending colimits to limits.
For example, using that extension, we get $X_\squareA\cong X_\triangle\times_{X_I}X_\triangle$.

\begin{theorem}\nonumber
Let $X_\bullet$ be a functor as above, and denote by the same letter its canonical extension to $\widehat \C$. 
Suppose that we are given isomorphisms
\[
\tau:X_\squareA\,\,\to X_\squareB
\qquad\text{and}\qquad \theta:X_\triangle\to X_\triangle
\]
making the following diagrams commute:
\[
\text{1)}\quad
\begin{matrix}
\xymatrix@C=.35cm{
X_\squareA
\ar[dr]_(.3){\stackrel{4}{\textstyle \Pi} \put(1,1){$\scriptstyle d'_i$}\,\,}\ar[rr]^\tau&&X_\squareB
\ar[dl]^(.3){\!\!\stackrel{4}{\textstyle \Pi} \put(1,1){$\scriptstyle d''_i$}}
\\
&X_I^4
}
\end{matrix}
\qquad
\parbox{6.5cm}{\small
where $d'_i:X_\squareA\to X_I$, $d''_i:X_\squareB\to X_I$, $i=1.\,.\, 4$, are induced by the four face inclusions $I\hookrightarrow\square$.
}\]
\[
\text{2)}\quad
\begin{matrix}
\xymatrix{
X_\squareA\ar[r]^\tau\ar[d]_{1\!/2}&X_\squareB\ar[d]^{1\!/2}\\
X_\squareA\ar[r]^\tau&X_\squareB
}
\end{matrix}
\qquad
\parbox{6.5cm}{\small
where $1\!/2:X_\squareA\to X_\squareA$ and $1\!/2:X_\squareB\to X_\squareB$
are induced by half turn rotation of the square.
}\]
\[
\text{3)}\quad
\begin{matrix}
\xymatrix{
X_\squareA\ar[r]^\tau\ar[d]_{1\!/4}&X_\squareB\ar[d]^{1\!/4}\\
X_\squareB&\ar[l]_\tau X_\squareA
}
\end{matrix}
\qquad
\parbox{6.5cm}{\small
where $1\!/4:X_\squareA\to X_\squareB$ and $1\!/4:X_\squareB\to X_\squareA$
are induced by rotation by a quarter turn.
}\]
\[
\text{4)}\hspace{.3cm}
\begin{matrix}
\xymatrix@C=.5cm{
\ar@{}[r]|{\displaystyle X_\pentagonB\phantom{.}}="a"&\ar@{}[r]|{\phantom{.}\displaystyle X_\pentagonC\phantom{\big)}}="b"
\ar "a";"b"^{\scriptscriptstyle \tau\times 1}&\\
X_\pentagonA\ar;"a"^{\scriptscriptstyle 1\times \tau} \ar[r]^{\scriptscriptstyle \tau\times 1} &
X_\pentagonE \ar[r]^{\scriptscriptstyle 1\times \tau} &
X_\pentagonD
\ar@{<-};"b"+<3pt,-9pt>_(.65){\scriptscriptstyle 1\times \tau}
}
\end{matrix}
\hspace{.3cm}
\parbox{6.5cm}{\small
note that ``$1\times \tau$'' and ``$\tau\times 1$'' only become well defined once we have axioms 1) and 2).}
\]
\[
\begin{array}{ll}
\text{5)}
&
d_i\circ\theta=d_{4-i}\vspace{.2cm}\\
\text{6)}
&
\theta\circ r=r^{-1}\circ\theta\vspace{.2cm}\\
\text{7)}
&
\theta^2=1,
\end{array}
\hspace{1cm}
\text{8)}
\quad
\begin{matrix}
\xymatrix@C=1.3cm{
X_\squareA\ar[r]^\tau\ar[d]_{\theta\times\theta}&X_\squareB\ar[d]^{\theta\times\theta}\\
X_\squareA\ar[r]^{\tau\circ 1\!/2}& X_\squareB
}
\end{matrix}
\hspace{2.6cm}
\]
then there is a natural action of $J'_n$ on the manifold that $X_\bullet$ associates to a triangulated $n$-gon.
(For example, one gets an action of $J'_8$ on $X\put(.2,-4){\psfig{file=triangulation.eps, width=0.27cm}}\hspace{.23cm}\phantom{.}$).
\end{theorem}

We now use the above theorem to equip $\A(n)_{>0}$ with a $J'_n$ action.
Indeed, the manifolds $X_\triangle:=\A(3)_{>0}$ and $X_I:=\A(2)_{>0}$ fit into a diagram (\ref{diagram}),
and so provide a functor $\widehat \C\to \{\text{manifolds}\}$.
The space associated to a triangulated $n$-gon is $\A(n)_{>0}$, as can be seen from the parameterization (\ref{bullets}).
We let $\tau$ be the composite 
\[
\tau:X_\squareA\,\,\stackrel{\scriptscriptstyle\sim}{\longrightarrow} \A(4)_{>0}\stackrel{\scriptscriptstyle\sim}{\longrightarrow} X_\squareB\,,
\]
and $\theta$ be the map sending $(F_1,F_2,F_3)\in(\Fl)^3$ to $(F_3^\bot,F_2^\bot,F_1^\bot)$, 
where the orthogonal of a flag $F$ is given by 
\(
(V_1,\ldots,V_{m-1})^\bot:=(V_{m-1}^\bot,\ldots,V_1^\bot),
\)
along with $\pm$ the obvious volume forms.
The axioms {\it 1)\,--\,8)} are then easy to check.

Both $\tau$ and $\theta$ are composites of cluster exchange relations.
But the action of $J'_n$ on $\A(n)_{>0}$ is not cluster
(it doesn't satisfy the Laurent phenomenon; it doesn't preserve the canonical presymplectic form).
The reason is that $\theta$ is actually the composite of a cluster map with an automorphism that negates the cluster matrix.
In particular, it negates the presymplectic form.

\end{talk}

\end{document}